\documentclass[12pt]{paper}

   \usepackage{cite}
 \usepackage{amsmath}  
\usepackage{amsfonts}
\usepackage{amssymb}
 
 \usepackage{centernot}
 \usepackage{epsfig}
\usepackage{graphics}
 \usepackage{color}
 
 \def\zHG0{H_{\Gamma_0}^1(\ZOMq)}
 \newcommand{\Dom}{{\rm dom}}
 \newcommand{\Ldue}{L^2\left (0,T;L^2(\Gamma_1)\right )}
  \newcommand{\ztil}{~}
\newcommand{\Mo}{\mathbb{M}}

\usepackage{color}
\usepackage{epsfig,cmmib57}
\usepackage{amssymb,amsmath,exscale}
 \numberwithin{equation}{section}

\newcommand{\ZSUno}{\sum _{n=1}^{+\ZIN}}

\newcommand{\ZOMq}{\Omega}
  
\def\zphin{\varphi_n(x)}
 
\newcommand{\ZA}{{\mathcal A}}

\newcommand{\zg}{\gamma}
 
 \newcommand{\sumZ}{\sum _{n\in \mathbb{Z}'}}
 
\newcommand{\intT}{\int_0^T}
\newcommand{\intt}{\int_0^t}
\newcommand{\ints}{\int_0^s}

\newcommand{\ZCD}{(\cdot)}
\newtheorem{Theorem}{Theorem}  
\newtheorem{Corollary}[Theorem]{Corollary}
\newtheorem{Lemma}[Theorem]{Lemma}

\newtheorem{Remark}[Theorem]{Remark}
\newtheorem{Definition}[Theorem]{Definition}

\newcommand{\zdiaform}{\mbox{~~\zdia}}

\newcommand{\zaa}{\alpha}
\newcommand{\ZEP}{\epsilon}

\newcommand{\zdia}{~~\rule{1mm}{2mm}\par\medskip}

\newcommand{\ZLA}{\label}
\newcommand{\ZIN}{\infty}
\newcommand{\zProof}{{\noindent\bf\underbar{Proof}.}\ }

\newcommand{\zzr}{{\rm I\hskip-2.1pt R}}

\newcommand{\ZD}{\;\mbox{\rm d}}

\author{
L. Pandolfi\thanks{Retired from the Dipartimento di Scienze Matematiche ``Giuseppe Luigi Lagrange'', Politecnico di Torino, Corso Duca degli Abruzzi 24, 10129 Torino, Italy (luciano.pandolfi@polito.it)}
}
\title{Controllability of a linear system with persistent memory via boundary traction\thanks{
This papers fits into the research program of the GNAMPA-INDAM and has been written in the framework of the   ``Groupement de Recherche en Contr\^ole des EDP entre la France et l'Italie (CONEDP-CNRS)''.}}
\begin{document}
 
 \maketitle 
 
 \begin{abstract}
We consider a linear viscoelastic system of Maxwell-Boltzmann type. Hence, viscosity contributes a memory term to the elastic equation. The system is controlled via the traction exerted on a   part $\Gamma_1$ of the boundary of the body. We prove that if the associated elastic system (i.e. the elastic system without memory) is exactly controllable then the viscoelastic system is exactly controllable too. This is similar to the known result when the boundary deformation is controlled, but the proof is far more delicate since controllability under boundary traction corresponds to the fact that  a certain sequence of functions is  a Riesz-Fisher sequence, but not a Riesz sequence.
 \end{abstract}
  
 \section{Introduction}
The deformation of a class of viscoelastic materials is described by the following equation\footnote{obtained using the MacCamy trick, described for example in~\cite{PandolfiLIBRO2014}, from the system $ w''=\Delta w + \intt N(t-s) \Delta w(s)\ZD s $.}
\begin{equation}\ZLA{eq:sistema}
w''=\Delta w +bw+\intt K(t-s) w(s)\ZD s
\end{equation}
with suitable initial and boundary conditions, as described below. In this equation, $w=w(x,t)$, $x\in\ZOMq$ (a bounded region of $\zzr^d$ with $C^2$ boundary), the apex denotes time derivatives, i.e.  $w''=w_{tt}$ and $\Delta $ is the Laplacian in the space variable $x$. Note that   we write $w=w(t)=w(x,t)$ as more convenient.

We use $\zg_0$ and $\zg_1$ to denote the traces on     $\Gamma=\partial \ZOMq$ of a function and, respectively, of its normal derivative.

Eq.~(\ref{eq:sistema}) is supplemented with initial and boundary   conditions:
\begin{equation}
\ZLA{iniBOUNDARYcond}
w(0)=w_0\,,\quad w'(0)=w_1\,,\qquad   \quad \left\{\begin{array}{l}
\zg_0 w=0 \quad {\rm on}\quad \Gamma_0\\
\zg_1w=f =\mbox{control on $ \Gamma_1 $} \,.
\end{array}\right.
\end{equation}

The \emph{associated wave equation} to Eq.~(\ref{eq:sistema}) is 
\begin{equation}
\ZLA{eq:associatedwave}
u''=\Delta u
\end{equation}
with the same initial and boundary conditions as $w$.

We shall be consistent in the use of $w$ and $u$ to denote  the solutions  respectively of the controlled equations~(\ref{eq:sistema})   and~(\ref{eq:associatedwave}). 
When needed, in order to stress the dependence on $f$ we write $w_f$ or $u_f$.
We shall also need to examine uncontrolled systems, i.e. $f=0$. In this case we use respectively $\psi$ and $\phi$ in the place of $w$ and $u$.

We are going to study controllability (in the space  discussed below) under the action of square integrable controls. More precisely, we are going to prove that controllability of the associated wave equation can be lifted to the system with memory. This fact is similar to the corresponding property when the control  acts in the Dirichlet boundary condition, but   the proof is   more delicate for a reason we shall see below.

The assumptions in this paper are as follows:
\begin{enumerate}
\item $\Gamma_0$ and $\Gamma_1$ are relatively open subsets of $\Gamma=\partial\ZOMq$ such that $\Gamma_0\cup\Gamma_1=\partial\ZOMq$ and $\overline \Gamma_0\cap\overline \Gamma_1=\emptyset$     in order to avoid the difficulties examined in~\cite{Grisvard1989}. 
We assume also $\Gamma_0\neq \emptyset$ (solely for the sake of simplicity, see Remark~\ref{remaSULLAfrontiera}).
 \item The memory kernel $K(t)$ is continuous.
\item we define:
\[
T_0=2\inf_{x_0\in\zzr^d}\left\{\sup_{x\in\ZOMq}\ |x-x_0|\right\}\,.
\]
\item   the part $\Gamma_1$ of the boundary is chosen in such a way that Theorems\ztil\ref{teo:controlASSOCIATEDwave} and\ztil\ref{teo:controlloDIRICHLET} below hold.
The existence of $\Gamma_1$ with this property has been proved in~\cite{LasieckaTriggianiTRACTIONcontrol1989,KOmornikLIBRO}, after the preliminary results   in~\cite{LIONSTOMO1-1988}.
We don't need to describe the geometric properties of $ \Gamma_1 $ which are used to prove controllability since the idea in this paper is as follows: the already estabilished property of controllability of the associated wave equation is inherited by the equation with memory.
  
\end{enumerate} 

Among the many results in~\cite{LasieckaTriggianiTRACTIONcontrol1989}, we single out the following one which deals with square integrable controls  (see also~\cite[Theorems~4.8 and~6.19]{KOmornikLIBRO} for the control time $T_0$):
\begin{Theorem}\ZLA{teo:controlASSOCIATEDwave}
Let
\[
 \zHG0 =\left \{\phi\in H^1(\ZOMq)\,:\quad \zg_0\phi=0\ {\rm on}\ \Gamma_0\right \}\,.
\]
When $\Gamma_1$ is a suitable part of $\partial\ZOMq$, for every $T>T_0$   and  for every $w_0$, $\xi  $ in $\zHG0$ and every $w_1$, $\eta$ in $L^2(\ZOMq)$ there exists a control $ f\in \Ldue$ such that
$u_f(T)=\xi$, $u'_f(T)=\eta$ ($u$ is the solution of~(\ref{eq:associatedwave}) with the initial/boundary conditions as in~(\ref{iniBOUNDARYcond})). 

\end{Theorem}

\begin{Remark}\ZLA{remaSULLAfrontiera}
Note that:
\begin{itemize}
\item  this result has to be properly interpreted, as explained below, since $(u_f(t), u'_f(t))$ does not evolve in $\zHG0\times L^2(\ZOMq)$ (unless ${\rm dim}\,\ZOMq=1$).  
 
\item if $\Gamma_0=\emptyset$ a similar result holds but the controllability spaces   are factor spaces, with respect to the spaces of the function with null integral mean. Our main result, 
Theorem~\ref{teo:controlMEMORY}, can be extended to the case $\Gamma_0=\emptyset$. The procedure is analogous to the one we used in the case $\Gamma_0\neq \emptyset$ and the extension is left to the reader.
\end{itemize}
\end{Remark}

We recall also the following result on controllability of the \emph{system with memory,} but \emph{controlled  via the deformation,} i.e. with control acting in the Dirichlet boundary condition (see~\cite{PandolfiLIBRO2014,PandolfiJINTEQAPPL2015,PandolfiPARMA2016}):
\begin{Theorem}\ZLA{teo:controlloDIRICHLET}
Let us consider Eq.~(\ref{eq:sistema}) but now the boundary condition is
\begin{equation}\ZLA{eq:boundaDiDirichletTeoINI}
\zg_0 w(x,t)=0 \ {\rm on}\ \Gamma_0\,,\qquad \zg_0 w(x,t)=f(x,t) \ {\rm on}\ \Gamma_1 \,.
\end{equation}
Let $T>T_0$ and let $\Gamma_1$ be as in Theorem~\ref{teo:controlASSOCIATEDwave}.  For every $w_0$, $\xi  $ in $L^2(\ZOMq)$ and every $w_1$, $\eta$ in $H^{-1}(\ZOMq)$ there exists a control $ f\in \Ldue$ such that
$w_f(T)=\xi$, $w'_f(T)=\eta$ (here $w$ solves~(\ref{eq:sistema}) with the initial   conditions  $w(0)=w_0$, $w'(0)=w_1$ and the boundary conditions~(\ref{eq:boundaDiDirichletTeoINI})).
 
\end{Theorem}
 
 \begin{Remark}
 We note:
 \begin{itemize}
 \item Theorem~\ref{teo:controlloDIRICHLET} holds in particular if $b=0$ and $K=0$, i.e. it holds for the associated wave equation~(\ref{eq:associatedwave}) controlled by the boundary deformation (see~\cite[Theorem~6.5]{KOmornikLIBRO}).
 \item we repeat that it is possible to choose $\Gamma_0$ and $\Gamma_1$ so that both the theorems~\ref{teo:controlASSOCIATEDwave}  and~\ref{teo:controlloDIRICHLET} hold.
 
 \item when the deformation, instead of the traction is controlled, as in Theorem~\ref{teo:controlloDIRICHLET},  the computation of $w_f(T)$ and $w_f'(T)$ is not a difficulty, since in this case $\left (w(t),w'(t)\right )\in C\left ([0,T];L^2(\ZOMq)\times H^{-1}(\ZOMq)\right )$.
  \end{itemize}
 
 \end{Remark}

The result we are going to prove is:

\begin{Theorem}\ZLA{teo:controlMEMORY}
Let $K $ be continuous, $T>T_0$ and let  $\Gamma_0$ and $\Gamma_1$ be such that both the theorems~\ref{teo:controlASSOCIATEDwave}  and~\ref{teo:controlloDIRICHLET} hold.  For every $w_0$, $\xi  $ in $\zHG0$ and every $w_1$, $\eta$ in $L^2(\ZOMq)$ there exists a control $ f\in \Ldue$ such that
$w_f(T)=\xi$, $w'_f(T)=\eta$.
\end{Theorem}
\begin{Remark}
As usual when proving controllability, we can assume null initial conditions: $w_0=0$, $w_1=0$. This will be done in this paper.\zdia
\end{Remark}
The definition of $w_f(T) $ and $w'_f(T) $ has to be explained    since, similar to the solution of the wave equation,  $(w_f(t), w'_f(t))$   does not evolve in $\zHG0\times L^2(\ZOMq)$.  This observation does not apply to the case ${\rm dim}\,\ZOMq=1$. In this case controllability has been studied in~\cite{NegrescuAPPLIED2016}.

The organization of the paper is as follows. We need to be very precise on the definition of the operators which are involved in the analysis of controllability of the associated wave equation. This is done in Sect.~\ref{sect:PreliWave}. The solutions and the corresponding operator for the system with memory are introduced in Sect.~\ref{sect:soluMemo} while controllability is proved in Sect.~\ref{sect:dimoControlMEMORY}. Notations are in Sect~\ref{sec:Notations}.
 
\subsection{\ZLA{sec:Notations}Notations and operators}
We introduce the following notation:
\[
H^\zaa_{\Gamma_0}(\ZOMq)=\left\{\begin{array}{l}
\left \{\phi\in H^\zaa(\ZOMq)\,:\ \zg_0\phi=0 \ {\rm on}\ \Gamma_0 \right \}\  {\rm if}\ \zaa>  1/2\\
H^\zaa(\ZOMq)\ {\rm if}\ \zaa\in(0,1/2)\\
\left (H^{-\zaa}_{\Gamma_0}(\ZOMq)\right )'\ {\rm if}\ \zaa<0 
\end{array}\right.
\]
 (the case $\zaa=1/2$ is not encountered in this paper).

We introduce the operator $A$ in $L^2(\ZOMq)$:
\begin{equation}
\ZLA{DefiOPERatoreA}
{\rm dom}\,A=\left \{ \phi\in H^2_{\Gamma_0}\,:\quad \zg_1\phi=0\ {\rm on}\ \Gamma_1\right \}\qquad A\phi=\Delta\phi 
\end{equation}
(note that the condition $\phi\in H^2_{\Gamma_0} $ does not impose conditions to the normal derivatives $\zg_1$ on $\Gamma_0$).
The operator $A$ is selfadjoint positive with compact resolvent
(regularity of $\partial\ZOMq$ is crucial for this property, see~\cite{JaksichSimon1992}) 
 and it is boundedly invertible since $\Gamma_0\neq\emptyset$. Let $\{\zphin\}$ be an orthonormal basis of $L^2(\ZOMq)$ whose elements are eigenvectors of $A$:
\[
A \zphin=-\mu_n^2\zphin\,.
\]
Note that the eigenvalues are not simple in general, but have finite multiplicity.

 We introduce
\[
\ZA=i\left (-A\right )^{1/2}\,,\qquad R_+(t)=\frac{e^{\ZA t}+e^{-\ZA t}}{2}\,,\qquad R_-(t)=\frac{e^{\ZA t}-e^{-\ZA t}}{2}\,.
\]
In fact, the operator $\ZA$ generates a $C_0$-group of operators.

It turns out that (see\cite[Sect.~2.1]{LasieckaTriggianiTRACTIONcontrol1989})
\[
\Dom\,\ZA=\zHG0\,.
\]

Finally we introduce the operator $G\in \mathcal{L}\left (L^2(\Gamma_1),L^2(\ZOMq)\right )$
\[
u=Gf\ \iff\ \left (\Delta u=0\ {\rm and}\  \left\{\begin{array}
{ll}
\zg_0 u=0& {\rm on}\ \Gamma_0\\
\zg_1 u=f& {\rm on}\ \Gamma_1
\end{array}\right.\right )
\]

It is known that $G$ takes values in $H^{3/2}_{\Gamma_0}(\ZOMq)\subseteq \Dom\,(-A)^{(3/4)-\epsilon}$ ($\epsilon>0$) which is compactly embedded in $L^2(\ZOMq)$.
In particular we have ${\rm im}\, G\subseteq \Dom\,\ZA$ (see~\cite[p.~195]{LasieckaTriggianiLIBROencVOL1}). Furthermore we note (see~\cite[Lemma~3.2]{LasieckaTriggianiTRACTIONcontrol1989}):
\begin{equation}\ZLA{eq:GsTelLaAsuBounda}
-G^*A\phi= \zg_0\phi _{|_{\Gamma_1}}      \quad \mbox{for every $\phi\in \Dom\, A$}\,. 
  \end{equation}

\section{\ZLA{sect:PreliWave}Preliminaries on the wave equation}
Here we report known properties on the wave equation with Neumann boundary conditions (see~\cite{LasieckaTriggianiCOSINE1981,LasieckaTriggianiNeumannREGULArityII1991}).
We consider the wave equation
\begin{equation}
\ZLA{eq:waveCOMPLETA}
u''=\Delta u+F\,,\qquad \left\{\begin{array}
{ll}
u(0)=u_0\,, & u'(0)=u_1\\
\zg_0 u=0& {\rm on}\ \Gamma_0\\
\zg_1u=f & {\rm on}\ \Gamma_1\,.
\end{array}\right.
\end{equation}
We assume $F\in L^2\left (0,T;L^2(\ZOMq)\right )$, $f\in \Ldue$ for every $T>0$. 
It is known that there exists $\zaa\in (0,1)$   such that for every $(u_0,u_1)\in H_{\Gamma_0}^\zaa(\ZOMq)\times   H^{ \alpha-1}_{\Gamma_0}(\ZOMq) $ problem~(\ref{eq:waveCOMPLETA}) admits a unique solution 
$u\in C\left ([0,T];H_{\Gamma_0}^\zaa(\ZOMq)\right )
\cap  
C^1\left ([0,T];H^{ \alpha-1}_{\Gamma_0}(\ZOMq) 
\right )$ 
and the transformation $(u_0,u_1,F,f)\mapsto u$ is continuous in the specified spaces. It is known that we can take $\zaa=(3/5)-\ZEP$ (any $\ZEP>0$); in particular $\zaa>1/2$ and so $1-\zaa<1/2$).   The values of $\zaa $ can be improved for special geometries 
  but in any case it will be $\zaa<1$, unless ${\rm dim}\,\ZOMq=1$
  (see~\cite{LasieckaTriggianiNeumannREGULArityII1991}   and~\cite[p.~739-740]{LasieckaTriggianiLIBROencVOL2}).

We need also an additional information on the special case $f=0$.

If $f=0$ then, as we stated already, the solution is denoted $\phi$ (instead of $u$) and it turns out that the map $(\phi_0,\phi_1,F)\mapsto u$ is linear and continuous from 
\[
\zHG0\times L^2(\ZOMq)\times L^2\left (0,T;L^2(\ZOMq)\right )\mapsto C\left ([0,T];\zHG0\right )\cap C^1\left ([0,T];L^2(\ZOMq)\right )\,.
\]

We shall use the following representation of the solutions, from~\cite{LasieckaTriggianiCOSINE1981}:
\begin{multline}\ZLA{Eq:rappreSOLUONDE}
u(t)=R_+(t)u_0+\ZA^{-1}R_-(t)u_1+\ZA^{-1}\intt R_-(t-s)F(s)\ZD s\\
-\ZA \intt R_-(t-s)G f(s)\ZD s
\end{multline}
and so
\begin{multline}\ZLA{Eq:rappreDERIVATAONDE}
u'(t)=\ZA R_-(t)u_0+ R_+(t)u_1+ \intt R_+(t-s)F(s)\ZD s\\
-A \intt R_+(t-s)G f(s)\ZD s\,.
\end{multline}
We repeat:
\begin{align*}
 {\rm im}\, G\subseteq H_{\Gamma_0}^{3/2 }(\ZOMq)\subseteq  \zHG0= \Dom\, \ZA \quad
  \mbox{so that}\quad   \ZA G\in\mathcal{L}\left (L^2(\Gamma_1),L^2(\ZOMq)\right )  \,.
\end{align*}
An integration by parts (justified in~\cite{PandolfiAMO2005})  shows:
\begin{Lemma}
Let $u_0=0$, $u_1=0$, $F=0$ and  $f\in C^1\left ([0,T];L^2(\Gamma_1)\right )$. 
Then $(u,u')\in 
C\left ([0,T];\zHG0\times L^2(\ZOMq)\right )$ 
and
\begin{align}
\nonumber u_f(t)&=Gf(t)-R_ +(t)Gf(0)-\intt R_+(t-s)Gf'(s)\ZD s\,,\\
\ZLA{eq:uCONfSmooth}u_f'(t)&=-\ZA R_-(t)Gf(0)-\ZA \intt R_-(t-s)Gf'(s)\ZD s 
\end{align}
so that $t\mapsto \left (u_f(t),u'_f(t)\right )\in C\left ([0,T];\zHG0\times L^2(\ZOMq)\right )$\,.
\end{Lemma}
 
 Hence, when $ f\in C^1 $ we can define both the maps 
\begin{equation}
\ZLA{eq:mappaLAMBDAdatireg}
\left\{\begin{array}{ll}
\Lambda\,: \ f\mapsto \left (u(T),u'(T)\right ) &\Ldue\to \zHG0\times L^2(\ZOMq)\,,\\
 \hat \Lambda\,:\ 
 f\mapsto \left (\ZA u(T),u'(T)\right )&\Ldue\to L^2(\ZOMq)\times L^2(\ZOMq)\,.
 \end{array}\right.
  \end{equation}
  These maps (with values in the spaces specified in~(\ref{eq:mappaLAMBDAdatireg}))  cannot be defined if $f$ is square integrable since in this case the function  $\left (u(t),u'(t)\right )$ evolves in a larger space.  
 We prove:
 
 \begin{Theorem}\ZLA{theo:closaElasticMAP} The maps $\Lambda$ and $\hat\Lambda$ on $ \Ldue $   to respectively $ \zHG0\times L^2(\ZOMq) $ and to $ L^2(\ZOMq)\times L^2(\ZOMq) $,
 originally defined when $ f $ is smooth,  are closable. 
\end{Theorem}
\zProof 
 It is sufficient that we prove closability of $  \Lambda$ since $\ZA$ is bounded and boundedly invertible from $\zHG0$ to $L^2(\ZOMq)$.
 
 Let $f_n\to 0$ in $\Ldue$, and $f_n\in C^1\left ([0,T];L^2(\Gamma_1)\right )$ so that $
 \left (u_{f_n}(T),u'_{f_n}(T)\right ) \in \zHG0 \times L^2(\ZOMq)
$   is well defined.

The sequence $\left  \{   \left ( u_{f_n}(T),u'_{f_n}(T)\right )\right \}$ in general does not converge in $\zHG0 \times L^2(\ZOMq) $. We must prove that \emph{if} it converges then it converges to $0$. We consider the sequence $\left  \{   u_{f_n}(T)\right \}$ of  the first components. The sequence of the velocities is treated analogously. 

Let   $  u_{f_n}(T)\to y\in \zHG0$ (the convergence is in the norm of $\zHG0$).
We noted already that $f\mapsto u_{f}(T)$ is continuous from $\Ldue$ to the 
 \emph{larger} space $H^{\zaa}_{\Gamma_0}(\ZOMq)$, and so $ u_{f_n}(T)\to 0$ in $H^{\zaa}_{\Gamma_0}(\ZOMq)$; The space $\zHG0$ is continuously embedded in $H^{\zaa}_{\Gamma_0}(\ZOMq)$ and so we see that $u_{f_n}(T)\to y$ in $H^{\zaa}_{\Gamma_0}(\ZOMq)$ too. And so it must be $y=0$.\zdia

 This result allows to extend  $\Lambda$ and $\hat \Lambda$ as closed operators to a certain dense subspace $\mathcal{F}$ of $\Ldue$:
$
f\in \mathcal{F}
$ when there exists a   sequence of smooth functions $f_n$ which converges to $f$ and such that $\left\{\Lambda f_n\right \}$  is convergent in $\zHG0\times L^2(\ZOMq)$. The limit is by definition $\Lambda f$ (the operator $\hat\Lambda$ is defined similarly). 

 It is clear that the vector $\Lambda f\in \zHG0$   does  not depend on the approximating sequence.

 \begin{Definition}{ 
From now on,  $\Lambda$ and $\hat\Lambda$ will denote these (minimal) closed extensions of the operators in\ztil(\ref{eq:mappaLAMBDAdatireg}), originally defined for smooth $f$.\zdia
 }
 \end{Definition}
 
 The result in Theorem~\ref{teo:controlASSOCIATEDwave} (reported from~\cite{LasieckaTriggianiTRACTIONcontrol1989,KOmornikLIBRO})  states that the operators $\Lambda$ and $\hat \Lambda$ (defined on $\mathcal{F}$) are surjective for every $T>T_0$.

We need the computation of the adjoints, and it is sufficient that we compute $\hat\Lambda^*$, which is closed and has dense domain, since $\hat\Lambda$ is closed.
So we can compute the adjoint in a dense subset of its domain, and then extend with the minimal closed extension.  Moreover, the computation of the adjoint can be done by restricting $\hat\Lambda$ to $C^1$ functions $f$ which are zero for $t=0$ and for $t=T$. Then, from~(\ref{eq:uCONfSmooth}), we have  
\[
u(T)= - \intT R_+(T-s)G f'(s)\,,\qquad u'(T)=- \ZA\intT R_-(T-s)Gf'(s)\ZD s\,.
\]

Let $ \xi$, $\eta$ belong to $\Dom\, A$ . Then   (see also~\cite[Lemma~3.3]{LasieckaTriggianiTRACTIONcontrol1989}):
\begin{equation}\ZLA{eq:preliAggiu}
\hat \Lambda^*(\xi,\eta) = - G^*A\underbrace{\left [R_+(T-s)\xi+\ZA^{-1}R_-(T-s)\left (\ZA \eta\right )\right ]}_{\phi(T-s)} 
\,.
\end{equation}
For example we compute 
\begin{multline}\ZLA{Eq:CompuAGGIuperLambdaSTAR}
\int_\ZOMq \ZA u(T)\eta\ZD x=- \int_\ZOMq\left [\ZA \intT R_+(T-s)Gf'(s)\ZD s\right ]\eta\ZD x \\
=- \intT\int_{\Gamma_1} f'(s)\left [G^*\ZA R_+(T-s)\eta\right ]\ZD\Gamma\ZD s \\
= - \intT\int_{\Gamma_1}  f(s)\left [G^*A  \left (\ZA^{-1}R_-(T-s) \ZA \eta\right )\right ]\ZD\Gamma\,\ZD s\,.
\end{multline}
Analogously
\begin{equation}
\ZLA{Eq:CompuAGGIuperLambdaSTARUno}
\int_{\ZOMq} u'(T)\xi\ZD x=-\intT\int_{\Gamma_1} f(s)G^*AR_+(T-s)\xi\ZD\Gamma\,\ZD s\,.
\end{equation}

Using~(\ref{eq:GsTelLaAsuBounda}), formula~(\ref{eq:preliAggiu}) is easily interpreted:   $\phi(t)$ 
in~(\ref{eq:preliAggiu}) solves
\begin{equation}\ZLA{EquadiPhi}
\phi''=\Delta\phi\qquad \phi(0)=\xi\,,\ \phi'(0)=\ZA\eta\,,\qquad\left\{\begin{array}{ll}
\zg_0\phi=0&{\rm on}\ \Gamma_0\\
\zg_1\phi=0&{\rm on}\ \Gamma_1\,. 
\end{array}\right.
\end{equation}
Hence, when $\xi$, $\eta$ belong to $\Dom\,A$ then $\hat\Lambda^*(\xi,\eta)=\zg_0\phi_{|_{\Gamma_1}}$.
And so
$
(\xi,\eta)\in\Dom\,\hat \Lambda^*
$
when there exists a sequence of smooth elements $(\xi^N,\eta^N) $
such that
\[\left\{
\begin{array}{l}
(\xi,\eta)=\lim (\xi^N,\eta^N) \ \mbox{in $L^2(\ZOMq)\times L^2(\ZOMq)$}\\
 \lim \hat \Lambda^*(\xi^N,\eta^N)=
\lim\left [ G^*A\phi^N\ZCD\right ]
\ \mbox{exists in $\Ldue$} 
\end{array}\right.
\]
(here $\phi^N$ solves~(\ref{EquadiPhi}) with data $\xi^N$ and $\eta^N$).

By definition $\hat\Lambda^*(\xi ,\eta )=\lim \hat \Lambda^*(\xi^N,\eta^N)$.
 
 Now we recall that $(\xi,\eta)\in \Dom\,\hat\Lambda^*$ when the function
 \[
f\mapsto \langle \hat \Lambda f,(\xi,\eta)\rangle _{L^2(\ZOMq)\times L^2(\ZOMq)} 
 \]
 is continuous on $\Ldue$. Looking at~(\ref{Eq:CompuAGGIuperLambdaSTAR}) and using ${\rm im }\,G\subseteq \Dom\,\ZA$ we see that
 
 \[
 \int_\ZOMq \eta \ZA u(T) \ZD x=\intT\int_\ZOMq\left (\ZA Gf \right )\ZA R_-(T-s)\eta\ZD x\,\ZD s
 \]
 is a continuous function of $f\in\Ldue$ when $\eta\in\Dom\,\ZA$. Treating~(\ref{Eq:CompuAGGIuperLambdaSTARUno}) analogously we get the first statement in the next lemma:

   \begin{Lemma}\ZLA{eq:Lemma:SULdominiodiLambda}
 Let  $ \xi\in\Dom\,\ZA $ and $ \eta\in\Dom\,\ZA $ then: 
\begin{enumerate}
\item 
  $ (\xi,\eta)\in\Dom\,\Lambda^*=\Dom\,\hat\Lambda^* $. 
   
\item   
The transformation
$(\xi,\eta)\mapsto \hat\Lambda^*(\xi,\eta)$  restricted to to $\zHG0\times \zHG0=\Dom\,\ZA\times\Dom\,\ZA$ is continuous.
 
  \end{enumerate}
   \end{Lemma}
  \zProof The first statement was noted already.  

The proof of the second statement is as follows:
Let  $\{(\xi^N,\eta^N)\}\in \Dom\,\ZA\times\Dom\,\ZA$  and let
\[
\|\xi^N-\xi\|_{\Dom\,\ZA}\to 0\,,\quad \|\eta^N-\eta\|_{\Dom\,\ZA}\to 0\,.
 \]
We must prove that
\[
\lim _{N\to+\ZIN}\hat\Lambda^*\left (\xi^N,\eta^N\right )= \hat\Lambda^*\left (\xi ,\eta \right ) \quad \mbox{in the norm of $\Ldue$}\,.
\]
Note  from~(\ref{eq:preliAggiu}):
\begin{multline*}
\hat\Lambda^*\left (\xi^N,\eta^N\right )=-G^*A\left [R_+(T-s)\xi^N+R_-(T-s)\eta^N\right ]\\
=-G^*\ZA
\left [R_+(T-s)\left (\ZA\xi^N\right )+R_-(T-s)\left (\ZA\eta^N\right )\right ]\,.
\end{multline*}
  The condition $\left  (\xi^N,\eta^N\right )\to(\xi,\eta)$ in $\Dom\,\ZA\times\Dom\,\ZA$ is the condition 
  $\left  (\ZA\xi^N,\ZA \eta^N\right )\to(\ZA\xi,\ZA\eta)$ in $L^2(\ZOMq)\times L^2(\ZOMq)$ (here we use $0\in\rho(\ZA)$, i.e. $\Gamma_0\not= \emptyset$) and $G^*\ZA$ is continuous on $L^2(\ZOMq)$. So,
  \begin{multline*}
\lim _{N\to+\ZIN}  \hat\Lambda^*\left (\xi^N,\eta^N\right )=-G^*A\left [R_+(T-s)\xi^N+R_-(T-s)\eta^N\right ]\\
=-\lim _{N\to+\ZIN}G^*\ZA
\left [R_+(T-s)\left (\ZA\xi^N\right )+R_-(T-s)\left (\ZA\eta^N\right )\right ]\\=
-G^*\ZA
\left [R_+(T-s) \ZA\xi +R_-(T-s) \ZA\eta^N \right ]=\hat\Lambda^*(\xi,\eta)\,.\zdiaform
\end{multline*}
 \begin{Remark}
 The second statement of the lemma can  be interpreted as follows: the map $(\xi,\eta)\mapsto \zg_0\phi(\cdot)$ (defined as a transformation on $\Dom\,A\times\Dom\,A$ with values in $\Ldue$) admits a continuous extension to $\Dom\,\ZA\times\Dom\,\ZA$.\zdia
 \end{Remark}

 \subsection{Fourier expansions}
 We shall need the expansion of $\hat\Lambda$ and $\hat \Lambda^*$  in series of the $\varphi_n$, the   orthonormal basis of $ L^2(\ZOMq) $ we already fixed,   whose elements are eigenfunctions of $A$.
 
In order to find an expansion of $\Lambda$ we write
\[
u(x,t)=\ZSUno \varphi_n(x) u_n(t)
\]
and it is easily seen that $u_n(t)$ solves
\[
u_n''=-\mu_n^2 u_n+\int_{\Gamma_1} \zg_0\varphi_n f\ZD\Gamma
\]
so that
\begin{align*}
u_n(t)&=\frac{1}{\mu_n}\intt\int_{\Gamma_1} \left (\zg_0\varphi_n \sin\mu_n s\right )f(x,T-s)\ZD\Gamma\,\ZD s\,,\\
u_n'(t)&=\intt\int_{\Gamma_1} \left (\zg_0\varphi_n \cos\mu_n s\right )f(x,T-s)\ZD\Gamma\,\ZD s\,.
\end{align*}

It follows that
 
\begin{multline}\ZLA{eq:PerESPAserieHATlambda}
\hat\Lambda f=\ZSUno\varphi_n(x) \left [ \left (
 \intT\int_{\Gamma_1} \left (\zg_0\varphi_n \sin\mu_n s \right )f(x,T-s)\ZD\Gamma\,\ZD s\right )\right .\,,\\
\left. \left (\intT\int_{\Gamma_1} \left (\zg_0\varphi_n \cos\mu_n s \right )f(x,T-s)\ZD\Gamma\,\ZD s\right ) \right]
\end{multline}
and the domain of $\hat\Lambda$ (i.e. also that of $\Lambda$) is the set of the functions $f\in\Ldue$ such that the elements in the bracket   constitute an $l^2$ sequence. This statement has to be precisely justified from the definition of the operators, and the justification is as follows. Let $ \{f^K\} $ be a sequence of smooth functions, $ f^K\to f\in\Dom\,\hat\Lambda $. Then ($\left \{\zaa_n\right \}$ and $\left \{\zaa_n^K\right \}$ are the brackets in~(\ref{eq:PerESPAserieHATlambda}), computed with $f$ and with $f^K$)
\[
\hat\Lambda f^k=\ZSUno \varphi_n(x)\zaa^K_n\to 
\ZSUno \varphi_n(x)\zaa _n=\hat\Lambda f\,.
 \]
Using the fact that $ \{\varphi_n\} $ is an orthonormal sequence, we see that $ \left \{ \zaa^K_n\right \}\in l^2 $ for every fixed $ K $ and $   \left \{ \zaa^K_n\right \}\to   \{ \zaa _n\}$
in $ l^2 $. So we have  $ \{\alpha_n\}\in l^2 $ i.e.
\begin{multline*}
\left \{  \left (
 \intT\int_{\Gamma_1} \left (\zg_0\varphi_n \sin\mu_n s \right )f(x,T-s)\ZD\Gamma\,\ZD s\right )\,, \right.\\
 \left. \left (\intT\int_{\Gamma_1} \left (\zg_0\varphi_n \cos\mu_n s \right )f(x,T-s)\ZD\Gamma\,\ZD s\right ) \right \}\in l^2 
 \end{multline*}
 (and conversely).
It is convenient to introduce the following operator~$\hat \Mo$: $\Ldue\to l^2$ 
\begin{multline*}
\hat\Mo f=\left \{\intT\int_{\Gamma_1} \left (\zg_0\varphi_n \cos\mu_n s \right )f(x,T-s)\ZD\Gamma\,\ZD s \right.\\
\left.+i\intT\int_{\Gamma_1} \left (\zg_0\varphi_n \sin\mu_n s\right )f(x,T-s)\ZD\Gamma\,\ZD  s\right \}\\
=\left \{
\intT\int_{\Gamma_1} \left (\zg_0\varphi_ne ^{i\mu_n s}\right ) f(x,T-s)\ZD\Gamma\,\ZD s 
\right \}\,.
\end{multline*}
The operator $\hat\Mo$, which is the \emph{moment operator} of the control problem for the wave equation, is surjective since the system is controllable. Unfortunately, it is not continuous and so we cannot conclude that the sequence $\left \{\zg_0\varphi_n e^{i \mu_n t}\right \}$ is a  Riesz sequence in $\Ldue$,  as for  the analogous sequence encountered when  the deformation (instead of the traction) is controlled. The term that it is used when the moment operator is surjective is that the sequence $\left \{\zg_0\varphi_n e^{i\mu_n t}\right \}$ is a \emph{Riesz-Fisher sequence} (see~\cite{Young}) and of course it is equivalent to the adjoint $\hat\Lambda^*$ being coercive. So now we expand the adjoint operator $\hat\Lambda^*$. Let (with $\{\xi_n\}\in l^2$, $ \{\eta_n\}\in l^2$)
\[
\xi=\ZSUno \xi_n\varphi_n(x)\,,\qquad \eta=\ZSUno \eta_n\varphi_n(x)\quad  \mbox{so that}\quad   \ZA\eta=\ZSUno \mu_n \eta_n\varphi_n(x) \,.
\]
The representation of the  solution $\phi$ of~(\ref{EquadiPhi}) is
\begin{multline}\ZLA{eq:repRESdrSolzPhi}
\phi(x,t)=\ZSUno\varphi_n(x)\left [\xi_n \cos\mu_n t+ \eta_n \sin\mu_n t\right ]\\
=\lim _{N} \sum _{n=1}^N \varphi_n(x)\left [\xi_n \cos\mu_n t+ \eta_n \sin\mu_n t\right ]\,.
\end{multline}
We observe that
\[
\xi^N=\sum _{n=1}^N\xi_n \varphi_n(x)\,,\quad \eta^N=\sum _{n=1}^N\eta_n \varphi_n(x)
\]
both belong to $\Dom\,A$ and can be used in the definition of $\hat\Lambda^*$.
So,  from~(\ref{eq:preliAggiu}) with $T-s$ replaced by $t$, 
\[ 
\hat \Lambda^*( \xi,\eta)  
=\lim _{N} \sum _{n=1}^N\left (- G^*A\varphi_n(x)\right )\left [\xi_n \cos\mu_n t+ \eta_n \sin\mu_n t\right ]\,:
\]
we have that $(\xi,\eta)\in\Dom\,\hat \Lambda^*$ when the limit exists in $\Ldue$ and then
\begin{multline}
\ZLA{eq:espadiLambdaSTAR}
\hat \Lambda^*( \xi,\eta)  
=\lim _{N} \sum _{n=1}^N\left (- G^*A\varphi_n(x)\right )\left [\xi_n \cos\mu_n t+ \eta_n \sin\mu_n t\right ]\\
= \sum _{n=1}^{+\ZIN}\left (- G^*A\varphi_n(x)\right )\left [\xi_n \cos\mu_n t+ \eta_n \sin\mu_n t\right ]
\,.
\end{multline}
      
Furthermore we note:
 
\begin{Lemma}
\ZLA{lemma:nuovoREGOLAR-ZA}
Let $(\xi,\eta)\in\Dom\, \hat\Lambda^*$. The series
\[
 \ZSUno \left (G^*A\varphi_n(x)\right )\left [
  \frac{\xi_n}{\mu_n} \sin\mu_n t -
   \frac{\eta_n}{\mu_n} \cos\mu_n t
 \right ]
\]
belongs to $H^1\left (0,T;L^2(\Gamma_1)\right )$ and the convergence of the series is in this space.
\end{Lemma}  
\zProof The convergence of the series is clear from Lemma~\ref{eq:Lemma:SULdominiodiLambda}, since the series correspond to $\left  (-\ZA^{-1}\eta,\ZA^{-1}\xi\right )\in\Dom\,\hat\Lambda^*$. The formal termwise computation of the derivative gives the series of $(\xi,\eta)$ which converges in $\Ldue$, since $(\xi,\eta)\in\hat\Lambda^*$ by assumption. So, the series belongs to $H^1\left (0,T;L^2(\Gamma_1)\right )$ and the partial sums converge in this space.\zdia

\section{\ZLA{sect:soluMemo}The solution of the system with persistent memory}

We define the solutions of the problem~(\ref{eq:sistema})-(\ref{iniBOUNDARYcond})  and of the corresponding problem for $\psi$, when $f=0$. In order to have a unified tratement, we assume that the initial conditions for $w$ are possibly non zero, as  in~(\ref{iniBOUNDARYcond}):    $w(0)=w_0$, $w'(0)=w_1$. Then, formally solving the wave equation~(\ref{eq:sistema}) perturbed by the affine term
\[
F(t)=bw(t)+\intt K(t-s)w(s)\ZD s
\]
we find
\begin{equation}\ZLA{eq:VolteDiw}\left\{\begin{array}{l}
\displaystyle  w( t)= u(t)+\ZA^{-1}\intt R_-(t-s)\left [bw(s)+\ints K(s-r)w(r)\ZD r\right ]\ZD s \\[2mm]
\displaystyle   w'(t)=u'(t)+\intt R_+(t-s)\left [bw(s)+\ints K(s-r)w(r)\ZD r\right]\ZD s
 \end{array}\right.
 \end{equation}
where
\[
u(t)= R_+(t)w_0+\ZA^{-1}R_-(t)w_1 -\ZA\intt R_-(t-s) Gf(s)\ZD s
\]
solves the associated wave equation with the same initial and boundary data.
 Note that the equation of $w'(t)$, i.e. the second line in~(\ref{eq:VolteDiw}), can also be written
\begin{multline}\ZLA{eq:VolteDiwDERIVATA}
w'(t)=u'(t)+
\ZA^{-1}\left [ R_-(t)bw_0+\intt R_-(t-s)K(s)w_0\ZD s\right ]\\
+\ZA^{-1}\intt R_-(t-s)\left [bw'(s)+\ints K(s-r)w'(r)\ZD r\right ]\ZD s\,.
\end{multline}

We noted that $u\in C\left ([0,T];H^\zaa_{\Gamma_0}(\ZOMq)\right )\cap
 C^1\left ([0,T]; H^{\zaa-1}_{\Gamma_0}(\ZOMq) \right )   
 $ for $\zaa>0$ small enough. 
So, from~\cite[p.~739-740]{LasieckaTriggianiLIBROencVOL2}, we have also
\[
  \left (u_f, u'_f\right )\in  C\left([0,T];\Dom\,(\ZA)^{\zaa }\times \left (\Dom\,(\ZA)^{1- \zaa }\right)'\right ) \,.
\]

 The Volterra integral operators in~(\ref{eq:VolteDiw}) leave these spaces invariant. So 
 we have also, for $\zaa\in(0,1)$ small enough,
\[
  \left (w , w' \right )\in  C\left([0,T];\Dom\,(\ZA)^{\zaa }\times \left (\Dom\,(\ZA)^{1- \zaa }\right)'\right ) 
\]
(continuous dependence on $w_0$, $w_1$ and $f$).

We repeat that  in order to get this property we use $\zaa>0$ small, in particular $ \zaa<1 $ and so $ \left (w(t),w'(t)\right ) $
\emph{ does not evolve in $\zHG0\times L^2(\ZOMq)  $.}

When $f=0$ we get the   solution $ \psi $  which   evolve in the same spaces as the solution $\phi$ of the associated wave equation according to the regularity of the initial conditions.

Now we define the operators $\Lambda_V$ and $\hat\Lambda_V$, which are analogous to the operators $\Lambda $ and $\hat\Lambda $.

When $f\in \mathcal{D}\left(\ZOMq\times (0,T)\right )$
the following definition makes sense:
\begin{align*}
\Lambda_V f&=\left (w(T),w'(T)\right )\in \zHG0\times L^2(\ZOMq)\,, \\
\hat \Lambda_V f&=\left (\ZA w(T),w'(T)\right )\in L^2(\ZOMq)\times L^2(\ZOMq)\,.
\end{align*}
As in Theorem~\ref{theo:closaElasticMAP}, we can see that these operators are closable and by definition their (minimal) closures are the operators $\Lambda_V$ and $\hat \Lambda_V$ used in the following  definition of controllability.

\begin{Definition}
Controllability of the system with memory is surjectivity of $\Lambda_V$ from $\Ldue$ to $\zHG0\times L^2(\ZOMq)$. Equivalently, the system with memory is controllable when    $\hat\Lambda_V$ is surjective from $\Ldue$ to $L^2(\ZOMq)\times L^2(\ZOMq)$.\zdia
\end{Definition}

We note:
\begin{Lemma}\ZLA{lemma:equalitDOMAINS}
The following properties hold:
\begin{enumerate}
\item\ZLA{lemma:equalitDOMAINSitem1}
the operators $\Lambda$ and $\Lambda_V$ have the same domain (and so  also  
$\hat \Lambda$ and $\hat \Lambda_V$ have the same domain).
 
\item\ZLA{lemma:equalitDOMAINSitem2}
The operators $  \left (\Lambda_V -\Lambda\right ) $
and $  \left (\hat \Lambda_V -\hat \Lambda\right ) $
 are compact from $L^2\left (0,T;L^2(\Gamma_1)\right )$ to, respectively,  $\zHG0\times L^2(\ZOMq)$ and to
 $L^2(\ZOMq)\times L^2(\ZOMq)$\,. 
 \item\ZLA{lemma:equalitDOMAINSitem3}
the operators $\Lambda^*$ and $\Lambda_V^*$ have the same domain (and so  also  
$\hat \Lambda^*$ and $\hat \Lambda_V^*$ have the same domain).
 \end{enumerate}
 \end{Lemma}
\zProof
We see from~(\ref{eq:VolteDiw}) and~(\ref{eq:VolteDiwDERIVATA}) that (here $w_0=0$, $w_1=0$)
\begin{equation}\ZLA{eqLambdaVcomePERTUR}
\Lambda_V f=\Lambda f+\mathcal{K}f
\end{equation}
where
\begin{multline*}
\mathcal{K}f=\left (\ZA^{-1}\intT R_-(T-s)\left [bw_f(s)+\ints K(s-r) w_f(r)\ZD r\right ]\ZD s,\right.\\
\left.  
\ZA^{-1}\intT R_-(T-s)\left [bw'_f(s)+\ints K(s-r) w'_f(r)\ZD r\right ]\ZD s \right)\,.
\end{multline*}

We noted that the transformation $f\mapsto \left (w_f,w'_f\right )$ is linear continuous from $\Ldue$ to $C\left ([0,T]; \Dom\,\ZA^\zaa\times \left (\Dom\,\ZA^{1-\zaa}\right )'\right )$ for e certain $\zaa\in (0,1)$. Hence $f\mapsto \mathcal{K}f$ is linear  \emph{ continuous and compact}   from $\Ldue$ to $\zHG0\times L^2(\ZOMq)$, since $ \ZA^{-1}$  is a compact operator.
The statements in the items~\ref{lemma:equalitDOMAINSitem1} and~\ref{lemma:equalitDOMAINSitem2} follow.

The statement in item~\ref{lemma:equalitDOMAINSitem3}  follows since $\Lambda_V^*=\Lambda_V+\mathcal{K}^*$ and $\mathcal{K}^*$ is continuous.\zdia

\begin{Remark}
{ Note that here we used compactness of the resolvent of $A$.\zdia}
\end{Remark}

It follows from Lemma~\ref{eq:Lemma:SULdominiodiLambda} that $ \Dom\,\Lambda_V^*=\Dom\,\hat\Lambda_V^*\supseteq \Dom\,\ZA\times\Dom\,\ZA  $.

Now we compute the adjoint  and its  expansions  in series of the eigenvectors $\{\varphi_n\}$.

In order to compute the adjoints we can again assume $f\in \mathcal{D}\left (\Gamma_1\times(0,T)\right )$ 
and $\xi$, $\eta$ smooth. Formulas~(\ref{eq:preliAggiu})  and~(\ref{EquadiPhi})
and the representation~(\ref{eqLambdaVcomePERTUR}) suggest that we consider 
\begin{equation}
\ZLA{EquadiPsi}
\psi''=\Delta\psi+b\psi+\intt K(t-s)\psi(s)\ZD s  
\end{equation}
with initial and boundary conditions
\begin{equation}\ZLA{eq:iniBOUNdperLambdaSTAR}
 \psi(0)=\xi\,,\ \psi'(0)=\ZA \eta\,,\qquad
\zg_0\psi=0\ {\rm on}\ \Gamma_0\,,\\ \zg_1\psi=0\ {\rm on}\ \Gamma_1\,.
\end{equation}
We assume that $ \xi $, $ \eta $   have finite expansions in series of the eigenfunctions $ \varphi_n $ and we compute
$ \hat \Lambda^*_V(\xi,\eta) $ in this case. Then we extend to the domain of   the minimal closure of the operator.

We multiply both the sides of~(\ref{eq:sistema}) with $\psi(T-t)$ and we integrate on $\ZOMq\times[0,T]$. 
We integrate by parts in time and space and (using $w(0)=0$, $w'(0)=0$) we get the equality: 
\begin{multline}\ZLA{eq:aggiDEB}
\intT\int_{\Gamma_1} \left (\zg_0 \psi(x,T-s)\right ) f(x,s)\ZD\Gamma\,\ZD s=\int_\ZOMq \xi w'( T)\ZD x+\int_\ZOMq \left (\ZA \eta\right )   w(T) \ZD x\\
=\langle
\left (\ZA w(T),w'(T)\right ),\left (\eta,\xi\right )
\rangle_{L^2(\ZOMq)\times L^2(\ZOMq)}\,.
\end{multline}

So, 
  
\begin{equation}\ZLA{eq:condiOrtogSec3}
\hat \Lambda^*_V(\xi,\eta)=\zg_0\psi(T-\cdot) =-G^*A \psi(T-\cdot)
\end{equation}
provided that $\psi$ solves~(\ref{EquadiPsi})  with   conditions~(\ref{eq:iniBOUNdperLambdaSTAR}) and 
$  ( \xi,\eta) \in L^2(\ZOMq)\times L^2(\ZOMq)$ are smooth, for example if they  have finite Fourier expansions.

We computed  $\hat\Lambda_V^*$    with smooth data but adjoint operators are    closed and   $\hat\Lambda_V^*$ is the   closed extension obtained as follows:
the elements of   $\Dom\,\hat \Lambda^*$  are those $(\xi,\eta)$ for which $\left \{-G^*A \psi^N\right \}$,  computed with smooth initial conditions $ (\xi^N, \eta^N )\to (\xi,\eta)$, is $\Ldue $-convergent and the limit is by definition  $ \hat\Lambda_V^*(\xi,\eta) $.

The computation of $ \Lambda^* $, defined on $ L^2(\ZOMq)\times \left (\zHG0\right )'$ is similar, but we don't need the details.

We repeat that as approximating sequences $\left\{(\xi^N,\eta^N)\right \}$ we can use sequences whose elements have finite expansions in series of the eigenfunctions $\varphi_n$, but the definition of the operators does not depend on the special sequence used.
 
 \begin{Remark}\ZLA{eq:RemaAGGiunSullExtRenD}
 The equality $\hat \Lambda^*_V(\xi,\eta)=\zg_0\psi(T-\cdot)$
  holds if $\xi$, $\eta$ have finite Fourier expansions. It holds also if they belong to $\Dom\,\ZA$ since in this case $\psi\in C\left ([0,T];\Dom\,\ZA\right )$ and, as we noted, $\zg_0$ is continuous on $\Dom\,\ZA$.\zdia 
 \end{Remark}

Finally we need the expansion of $\psi$ in series of the eigenfunctions $\varphi_n$. We consider the solution of~(\ref{EquadiPsi}) with conditions~(\ref{eq:iniBOUNdperLambdaSTAR}).
As $\xi\in L^2(\ZOMq)$, $\ZA\eta\in   \left(\Dom\,\ZA\right )' $ we have
\[
\xi(x)=\ZSUno\xi_n\varphi_n(x)\,,\qquad \ZA\eta(x)=\ZSUno \mu_n\eta_n \varphi_n(x)\,,\qquad \{\xi_n\}\,,\ \{\eta_n\}\in l^2\,.
\]
Hence:
\begin{align*}
  &\psi(t)=\ZSUno\varphi_n(x)\psi_n(t)\\
&\psi_n''= -\mu_n^2\psi_n+b\psi_n+\intt K(t-s)\psi_n(s)\ZD s\\
&\psi_n(0)=\xi_n\,,\qquad \psi_n'(0)=\mu_n\eta_n\,.
\end{align*}
 
So we have the following Volterra integral equation for $\psi_n(t)$:
\begin{multline}\ZLA{eq:soluPsi}
\psi_n(t)=\xi_n\cos\mu_n t+\eta_n\sin\mu_n(t)\\+\frac{1}{\mu_n} \intt \left [b\psi_n(s)+\ints K(s-r)\psi_n(r)\ZD r\right ]\sin\mu_n (t-s)\ZD s\,.
\end{multline}

So, we have the following equality if $ (\xi,\eta)\in\Dom\,\hat\Lambda_V^* $ (we replace $T-s$ with $t$);
 \begin{multline}\ZLA{eqExpeSeriHAtLaV}
 \hat\Lambda_V^*(\xi,\eta)=\lim _{N}\left [-G^*A\left (\sum _{n=1 }^N\varphi_n(x)\psi_n(t)\right )\right ]\\
 = \lim _{N}\sum _{n=1}^N
\left (-G^*A \varphi_n(x)\right )\psi_n(t)
 =\ZSUno \left (-G^*A \varphi_n(x)\right )\psi_n(t)
 \end{multline}
 (convergence in $\Ldue$).

Finally, also the analogous of the last statement in lemma~\ref{eq:Lemma:SULdominiodiLambda}  holds, with analogous proof: 
\begin{Lemma}\ZLA{lemmaDellaesteTRACCIAvisco}
The map $(\xi,\eta)\mapsto \hat\Lambda^*(\xi,\eta)$ restricted  to $\Dom\,\ZA\times \Dom\,\ZA=\zHG0\times\zHG0$ (as a map with value in $\Ldue$) is continuous. 
\end{Lemma}

 \begin{Remark}\ZLA{eq:REMAlimitPsi}
 {
 Gronwall inequality applied to~(\ref{eq:soluPsi}) shows that for every $ T>0 $ there exists $ M=M_T $ such that for every $ t\in[0,T] $ and every $n$ we have
 \[ 
 |\psi_n(t)|\leq M\left (\|\xi\|_{L^2(\ZOMq)}+\|\eta\|_{L^2(\ZOMq)}\right ). 
  \] 
  The number $M$ does not depend on $n$.\zdia
 }
 \end{Remark}
  
\section{\ZLA{sect:dimoControlMEMORY}The proof that the system with memory is controllable}

Let use put 
\[
R_V={\rm im}\,\Lambda_V\,,\qquad \hat R_V={\rm im}\,\hat\Lambda_V
\]
The fact that $\Lambda$ is surjective and $\Lambda-\Lambda_V$ is compact implies
\begin{Lemma} 
$R_V$ and $\hat R_V$ are closed with finite codimension (respectively in $\zHG0\times L^2(\ZOMq)$  and in $L^2(\ZOMq)\times L^2(\ZOMq)$).
\end{Lemma}
Hence, in order to prove controllability  it is sufficient to prove \emph{approximate controllability} i.e.  that the subspace $R_V$ is dense in $\zHG0\times L^2(\ZOMq)$, or that $\hat R_V$ is dense in $L^2(\ZOMq)\times L^2(\ZOMq)$. We   prove 
$\left  [\hat R_V\right ]^\perp =0$, i.e. we prove that 
 if $\hat \Lambda_V^* (\xi,\eta)=0$  then $(\xi,\eta)=0$.
 
Using~(\ref{eqExpeSeriHAtLaV}) we see that $\hat \Lambda_V^* (\xi,\eta)=0$ is the condition

\begin{equation}\ZLA{eq38RiscriINIpag14}
 \ZSUno \left (G^*A \varphi_n(x)\right )\psi_n(t) 
=0  \quad \mbox{(convergence  in $\Ldue$)}\,.
\end{equation}

 Our goal is the proof that condition~(\ref{eq38RiscriINIpag14}) implies $\xi=0$, $\eta=0$.

    The proof relies 
    on the following corollary to Theorem~\ref{teo:controlloDIRICHLET}. Note that in this corollary the space $H^1_0(\ZOMq)$, and not $\zHG0$, is used:  
    \begin{Corollary}\ZLA{COROateo:controlloDIRICHLET}Let $T>T_0$ and let $\Gamma_0$ and $\Gamma_1$ be as in Theorem~\ref{teo:controlloDIRICHLET}.
    Let $\psi$ solve~(\ref{EquadiPsi}) with conditions
    \[ 
    \psi(0)=\psi_0\in H^1_0(\ZOMq)\,,\quad \psi'(0)=\psi_1\in L^2(\ZOMq)\,,\qquad \left\{
    \begin{array}
    {ll}
    \zg_0\psi=0 &\mbox{on $\Gamma=\partial\ZOMq $}\\
    \zg_1\psi=0 &\mbox{on $\Gamma_1$.}
    \end{array}
    \right.
    \]
    Then $\psi(t)=0$ and so also $\psi_0=0$, $\psi_1=0$.
    \end{Corollary}
    The proof is in~\cite{PandolfiLIBRO2014,PandolfiJINTEQAPPL2015,PandolfiPARMA2016}.
\begin{Remark}\ZLA{rema:annHilatDirichl}
The assumptions on $\psi$ in Corollary~\ref{COROateo:controlloDIRICHLET} is the condition that $(\xi,\eta)$ annihilates the reachable set (in $L^2(\ZOMq)\times H^{-1}(\ZOMq)$) when the square integrable control acts on the deformation (i.e. in the Dirichlet boundary condition). It is also true that  (when $T>T_0$)   the converse implication holds, thanks to a compactness/unicity argument, but we are not going to use the converse implication.\zdia
\end{Remark}

Furthermore, we shall use the following result,  whose proof is postponed:
 \begin{Lemma}\ZLA{lemma:UlteREGOortog} Let $T>T_0$. If $(\xi,\eta)\in L^2(\ZOMq)\times L^2(\ZOMq)$ and  $(\xi,\eta)\perp \hat R_V$ then $(\xi,\eta)\in \Dom\,\ZA\times \Dom\,\ZA$.
 \end{Lemma}
 
 Granted this result, it is easy to see that $\xi=0$ and $\eta=0$ if $(\xi,\eta)\perp \hat R_V$. In fact, 
 Eq.~(\ref{EquadiPsi}) has now initial conditions $\psi(0)=\xi\in H^1_0(\ZOMq)$, $\psi'(0)=\ZA\eta\in L^2(\ZOMq)$. Hence $\psi$ evolves in $\zHG0$ and satisfies the following boundary conditions:
  
 \begin{equation}\ZLA{eq:ListaCondiBORDO}
\left\{\begin{array} {ll}
 \left\{\begin{array}{ll}
 \zg_0\psi=0\ {\rm on}\ \Gamma_0&\mbox{from~(\ref{eq:iniBOUNdperLambdaSTAR})}\\
 \zg_0\psi=0 \ {\rm on}\ \Gamma_1& \mbox{orthogonality condition}
 \end{array}\right.\\
 \zg_1\psi=0\ {\rm on}\ \Gamma_0\qquad \mbox{from~(\ref{eq:iniBOUNdperLambdaSTAR}).}
\end{array}\right.
 \end{equation}
 These properties are the conditions that $(\xi,\ZA\eta)$   annihilates the reachable set in $L^2(\ZOMq)\times H^{-1}(\ZOMq)$ of the control system~(\ref{eq:sistema}) \emph{with square integrable control in the Dirichlet boundary condition,} see Remark~\ref{rema:annHilatDirichl}.
 
Theorem~\ref{teo:controlloDIRICHLET} implies (via Corollary~\ref{COROateo:controlloDIRICHLET}) $\psi_0=\xi=0$, $\psi_1=\ZA\eta=0$ and so $\left [\hat R_V\right ]^\perp=0$, as we wished to prove. 

In fact, there are two    points to clarify:
\begin{itemize}

\item if $(\xi,\ZA\eta)$   annihilates the reachable set in $L^2(\ZOMq)\times H^{-1}(\ZOMq) $ \emph{of the system controlled via the Dirichlet boundary condition} then it must be $\xi\in H^1_0(\ZOMq)$, i.e. it must be $\zg_0\xi=0$ on the entire $\Gamma=\partial\ZOMq$.
Instead, we know from Lemma~\ref{lemma:UlteREGOortog} that $\xi\in\Dom\,\ZA=\zHG0$. The property $\zg_0\xi=0$ on the entire boundary of $\ZOMq$ follows since $\psi(t)\to\psi(0)=\xi$ in $\zHG0$, hence in the norm of $H^1(\ZOMq)$. We use again continuity of the trace   $\zg_0$    from $H^1(\ZOMq)$ to $L^2(\partial\ZOMq)$ and $\zg_0\psi(t)=0$ on $\partial\ZOMq$ (from~(\ref{eq:ListaCondiBORDO})). Passing to the limit we get $0=\zg_0\psi(0)=\zg_0\xi$ on $\partial\ZOMq$.
\item the orthogonality condition~(\ref{eq38RiscriINIpag14}) has been written $\left (\zg_0\psi\right )_{|_{\Gamma_1}}=0$   thanks to Lemma~\ref{lemma:UlteREGOortog} and Remark~\ref{eq:RemaAGGiunSullExtRenD}.
\end{itemize}

\emph{In conclusion, in order to complete the proof of Theorem~\ref{teo:controlMEMORY} we must prove Lemma~\ref{lemma:UlteREGOortog}.} The proof relies on the following result, whose proof is similar to the proof of Lemma~3.4 in~\cite{PandolfiLIBRO2014}. It is reported for completeness.
 
\begin{Lemma}\ZLA{eq:LEmmaDAllaAPPEnd}
Let $K$ be a Hilbert space and let $\{\mu_n\}$ be a sequence of real numbers. Assume that $\left \{k_n e^{i\mu_n t}\right \}$ is a Riesz-Fisher sequence in $L^2(0,T;K)$ and that $\{\zaa_n\}\in l^2$ is a sequence of complex numbers such that
\[
H(t)=\sum \zaa_n k_ne^{i\mu_n t}\in H^1\left ([0,T+h];K\right )\,,\qquad h>0\,.
\]
Then, $\{\mu_n \zaa_n\}\in l^2$.
\end{Lemma}
\zProof
We know from\ztil\cite[Proposition~IX.3]{Brezis}: let  $H\in H^1(0,T+h_0;K)$ and $0<h<h_0$ then there exists $C=C(H)>0$ independent of $h$ such that
\begin{equation}\ZLA{eq:Chap2:per:Probem:Chapter2:RappiNCRE}
\left |
\ZSUno  \zaa_n\mu_n   \frac{e^{i\mu_n  h}-1}{\mu_n  h}e^{i \mu_n  t}k_n
\right |^2_{L^2(0,T;K)}=     \left | \frac{H(t+h)-H(t)}{h}\right |^2_{L^2(0,T;K)}\leq C\;.
\end{equation}
The proof in this reference is for real valued functions, but it is easily adapted to Hilbert valued functions.

 Using the fact that $ \left \{e^{i\mu_n  t}k_n\right \} $ is a Riesz-Fisher sequence in $L^2(0,T;K)$  we see that

  \begin{multline*}
 \ZSUno  \left |
 \zaa_n\mu_n  \frac{e^{i\mu_n  h}-1}{ \mu_n  h}
 \right |^2\leq
 \frac{1}{m_0}\left |
\ZSUno  \zaa_n\mu_n   \frac{e^{i\mu_n  h}-1}{\mu_n  h}e^{i \mu_n  t}k_n
\right |^2_{L^2(0,T;K)}
 \\
 = \frac{1}{m_0}\left  | \frac{H( t+h)-H( t)}{h}\right  |^2_{L^2(0,T;K)}
\leq C/m_0\,.
  \end{multline*}
 The last equality holds for $ h $ ``small'', say  if $|h|<h_0/2$.

Let $ s $ be real. There exists $ s_0>0 $ such that:
 \[
 \left |  \frac{e^{i s}-1}{s}\right | ^2=\left ( \frac{\cos s-1}{s} \right )^2+\left ( \frac{\sin s}{s}  \right )^2>\frac{1}{2}\quad \mbox{for $ 0<s<s_0 $.}
  \]
Then we have, for every $ h\in (0,h_0/2)$,

\[
 \sum _{  \mu_n <s_0/h   }|\zaa_n\mu_n |^2\leq
2\ZSUno  \left |
 \zaa_n\mu_n   \frac{e^{i\mu_n  h}-1}{ \mu_n  h}
 \right |^2\leq2\frac{C}{m_0} \,.
\]
So,  $\{\zaa_n\mu_n \}\in l^2$ as wanted.\zdia

\paragraph{\underline{The proof of Lemma~\ref{lemma:UlteREGOortog} and so of Theorem~\ref{teo:controlMEMORY}}}

We introduce the following notations:
 \begin{align*}
& f^{(*0)}*g  =g \,,\quad
f^{(*1)}*g= f*g=\intt f(t-s)g(s)\ZD s\,,\\ 
&f^{(*n)}*g=f*\left (f^{(*(n-1))}*g\right ) \quad k_n=G^*A\phi_n\in L^2(\Gamma_1)\,,\\
&S_n=\sin\mu_n t\,,\qquad C_n=\cos\mu_n t\,,\qquad  E_n=e^{i\mu_n t}\,.
\end{align*}
The right hand side of~(\ref{eq:soluPsi}) is
 \begin{equation}\ZLA{eq38RiscrittaCONVOLante}
\psi_n=\underbrace{\xi_nC_n+\eta_n S_n}_{\phi_n}+\frac{1}{\mu_n}   \left (bS_n+K*S_n\right ) *\psi_n \,.
\end{equation}
 \begin{Remark}[On the notations]
{ 
We use $\varphi_n=\varphi_n(x)$ to denotes the eigenfunctions of $A$ while
$\phi_n=\phi_n(t)$ denotes the function in~(\ref{eq38RiscrittaCONVOLante}), which is the $n$-th component  of the solutions $\phi=\phi(x,t)$ of~(\ref{EquadiPhi}),
\[
\phi(x,t) 
=\ZSUno \varphi_n(x)\phi_n(t) 
\] 
(see~(\ref{eq:repRESdrSolzPhi})).
Note also that for simplicity we shall write $k_n$ in the place of $ G^*A\varphi_n$, since we shall use lemma~\ref{eq:LEmmaDAllaAPPEnd}.\zdia
}
\end{Remark}
 
  It is convenient to rewrite~(\ref{eq38RiscrittaCONVOLante}) in the form
  \begin{equation}
\ZLA{eq38RiscrittaCONVOL}
  \psi_n=\underbrace{c_n E_n+\bar c_n E_{-n}}_{\phi_n}+\frac{1}{\mu_n} \underbrace{ \left (bS_n+K*S_n\right )}_{G_n}*\psi_n\,\quad c_n=\frac{\xi_n-i\eta_n}{2}\,.
  \end{equation}
   $N$ steps of the Picard iteration give the following formula for $\psi_n(t)$:
 \begin{equation}
 \ZLA{eq:dopoNpassi}
 \psi_n=\phi_n+\frac{1}{\mu_n}G_n*\phi_n+\cdots+\frac{1}{\mu_n^N}G_n^{(*N)}*\phi_n+\frac{1}{\mu_n^{N+1}}G_n^{*(N+1)}*\psi_n\,.
 \end{equation}

We introduce the notations
\[
\mathbb{Z}'=\mathbb{Z}\setminus\{0\}\,,\qquad 
\mu_{-n}=-\mu_n  \qquad k_{-n}=k_n
\,,\qquad c_{-n}=\bar c_n\,.
\]
In order to prove Lemma~\ref{lemma:UlteREGOortog} we must prove that $c_n=\tilde c_n/\mu_n$, $\left \{\tilde c_n\right \}\in l^2\left (\mathbb{Z}'\right )$.

Using~(\ref{eq:dopoNpassi}), the condition of orthogonality~(\ref{eq38RiscriINIpag14}) can be written as follows:

\begin{multline}\ZLA{eq:SerOrthoDistrib}
 \sumZ k_n E_n(t) c_n+\ZSUno k_n\left [\sum _{k=1}^N\frac{1}{\mu_n^k}\left [bS_n+K*S_n\right]^{(*k)}*\phi_n\right ]\\
 =-\ZSUno k_n\frac{1}{\mu_n^{N+1}}\left [bS_n+K*S_n\right ]^{*(N+1)}  *\Psi_n\,.
\end{multline}
  The reasons why it is correct to distribute the series as above, provided that $N$ is large enough, are as follows:
 
  \begin{itemize}
    \item the series~(\ref{eq38RiscriINIpag14})  converges in $\Ldue$  since $(\xi,\eta)\in \Dom\, \Lambda_V^*$;
  \item the   series on the right side of~(\ref{eq:SerOrthoDistrib})  converges if $N$ is sufficiently large, since:
\begin{itemize}
\item the sequence $\{\psi_n(t)\}$ is bounded on $[0,T]$, see Remark~\ref{eq:REMAlimitPsi}.
\item there exists a contant $C$, which depends on $\ZOMq$ such that (see~\cite{HasselBARNEttPROCSYMP2012})
\begin{equation}
\ZLA{stimaTRACCEautovettori}
\| k_n\|_{L^2(\Gamma)}=
\|\zg_0\varphi_n\|_{L^2(\Gamma)}\leq C\sqrt[3]{\mu_n}  \,.
\end{equation}
\item $\mu_n>c n^{1/d}$ where $d= {\rm dim}\,\ZOMq$, see~\cite{JaksichSimon1992} and note that we denoted $-\mu_n^2$ the eigenvalues. Note that the (piecewise) regularity of $\partial\ZOMq$ is crucial for this estimate (see~\cite{JaksichSimon1992}).
\end{itemize}

  Thanks to this property, this series even converges to a $C^1$ function, provided that $N$ is large enough.
  \item 
 The first series on the left of~(\ref{eq:SerOrthoDistrib}) converges in $\Ldue$,  since $(\xi,\eta)\in \Dom\, \Lambda_V^*=\Dom\, \Lambda^*$.
   \end{itemize}
   
  \emph{From now on, the number $N$ of the steps of the Picard iteration is fixed, so large that the series on the right side of~(\ref{eq:SerOrthoDistrib}) converges to a $C^1$ function.} 
  We prove that the intermediate series can be distributed on its addenda, and converges to an $H^1$ function. The critical case is the case $k=1$.
    Using 
  \[
S_n*C_n=\frac{t}{2}S_n\,,\qquad S_n*S_n=\frac{1}{2\mu_n}S_n-\frac{t}{2}C_n  
  \] 
it is easily seen that
  \begin{multline*}
\ZSUno k_n \frac{1}{\mu_n}\left [bS_n+K*S_n\right ]*\phi_n\\
=b\left [\underbrace{\frac{t}{2}\left (
\ZSUno k_n\left [
-C_n\frac{\eta_n}{\mu_n}+S_n\frac{\xi_n}{\mu_n} 
\right ]
\right )}_{\fbox{1}} +\underbrace{\ZSUno k_nS_n\frac{\eta_n}{2\mu_n^2}}_{\fbox{2}}\right ]\\
+K*\left (
\fbox{1}+\fbox{2}
\right )\,.
  \end{multline*}
 In fact, the series in $\fbox{1}$  converges since $ \ZSUno \frac{\xi_n}{\mu_n}\varphi_n$ and $ \ZSUno\frac{\eta_n}{\mu_n}\varphi_n$ belong to $\Dom\,\ZA$, hence to $\Dom\,\Lambda^*$ (this is the first statement in 
 lemma~\ref{eq:Lemma:SULdominiodiLambda}). 
 So, from Lemma~\ref{lemma:nuovoREGOLAR-ZA},
it converges to an $H^1(0,T;L^2(\Gamma))$ function because  we are using  $(\xi,\eta)\in \Dom\, \hat\Lambda^*$.   For a stronger reason also the   series \fbox{2} converges to an $H^1$ function too. In fact,   $\{\eta_n/\mu_n^2\}$ are the Fourier coefficients of an element in $\Dom\,\ZA^2=\Dom\,A$. And so, the last term, which is the convolution of $K$ with an $H^1$-function, is of class $H^1$ too.

 The terms with $k\geq 2$ are treated analogously. 
 
 Hence, $ \sumZ k_n E_n(t) c_n\in H^1\left (0,T;L^2(\Gamma_1)\right )$ and we know that $\left  \{k_n E_n(t)\right \}$ is a Riesz-Fisher sequence in this space. Hence,
 \[
c_n=\frac{\tilde c_n}{\mu_n}\,,\qquad \left \{\tilde c_n\right \}\in l^2\left (\mathbb{Z}'\right )
 \]
 from Lemma~\ref{eq:LEmmaDAllaAPPEnd}. This is the result we wanted to achieve, see the statement of  Lemma~\ref{lemma:UlteREGOortog}, and completes the proof of controllability.

   \bibliography{bibliomemoria}{ }
  \bibliographystyle{plain}
 \enddocument